\documentclass[11pt, leqno]{article}
\usepackage{a4, indentfirst, latexsym,amssymb}
\usepackage{t1enc}
\usepackage{color}

  \newtheorem{Th}{Theorem}[section]
  
  \newtheorem{Prop}[Th]{Proposition}

  \newtheorem{Lem}[Th]{Lemma}

  \newtheorem{Rem}[Th]{Remark}

  \def\d{\,\mathrm{d}}

  \def\part{\partial}

  \def\R{\mathbb{R}}

  \newcommand{\lma}{\lambda}
  
  \newcommand{\eps}{\varepsilon}
  \newcommand{\be}{\begin{equation}}
  \newcommand{\ee}{\end{equation}}

\begin{document}

\begin{center}

{\huge  Averaging principle and hyperbolic evolution
equations}

Aleksander \'{C}wiszewski \footnote{
The research supported by the MNiSzW Grant no. N N201 395137.}\\

{\em Faculty of Mathematics and Computer Science\\
Nicolaus Copernicus University\\
ul. Chopina 12/18, 87-100 Toru\'n, Poland }\\
e-mail: aleks@mat.umk.pl\\

Version  12.04.2011

\end{center}

\begin{abstract}
An averaging principle is derived for the abstract nonlinear
evolution equation where the almost periodic right hand-side is a
continuous perturbation of the time-dependent family of linear
operators determining a linear evolution system. It generalizes classical Henry's results for perturbations of sectorial operators
on fractional spaces.
It is also proved that the main hypothesis of the nonlinear averaging principle is
satisfied for general hyperbolic evolution equations introduced by Kato.
\end{abstract}

\section{Introduction}
We are concerned with the limit behavior with regard to $\lma\to
0^+$ of evolution systems of the form
$$
\dot u(t)=  A(t/\lma) u(t) + F(t/\lma, u(t)), \, \, t>0,
\leqno{(P_\lma)}
$$
where $\{ A(t)\}_{t\geq 0}$ is a family of operators generating
$C_0$ semigroups of bounded linear operators on a Banach space
$E$, $F:[0,+\infty) \times E \to E$ is a continuous map satisfying
the local Lipschitz condition with respect to the second variable
and $\lambda>0$ is a parameter. The so-called
{\em averaging principle} is a well known tool in the theory of
ordinary differential equations, i.e., when $E$ is finite
dimensional. Roughly speaking, it says that if $F$ is periodic in
time, then trajectories of $\dot u (t)=F(t/\lma,u(t))$ converge to
trajectories of the averaged equation as $\lma\to 0^+$ (see
\cite{Bogoliubov-Mitropolsky}). This averaging idea is of
importance when studying qualitative behavior of nonautonomous
equations. It enables to perceive the dynamics of a nonautonomous
equation in terms of the related averaged one. For instance, by
this approach, one may examine global attractors for dissipative
equations, periodic solutions and other dynamic features such as
bounded or recurrent solutions. Therefore extending the method to
infinite dimension and applying it to partial differential
equations is a natural and vital issue attracting much attention.
The averaging principle in the infinite dimensional case was obtained by
Henry \cite{Henry} who assumed that the (independent of time) operator $A$ is a sectorial one on a Banach space $E$ and
$F:[0,+\infty)\times E^{\alpha} \to E$, where $E^{\alpha}$, $0\leq \alpha <1$, is the fractional power space determined by $A$,
is bouneded and continuous.
Averaging for time dependent (set-valued) perturbations of a $C_0$ group generator was considered by Kamenskii,
Obukhovskii and Zecca in \cite{Ka-Obu-Zec}, where $A$ was a $C_0$ semigroup generator and $F$ was an upper semicontinuous $k$-set conctraction
with respect to a measure of noncompactness.
Averaging principle, in the context of attractors and Conley-Rybakowski index, for
parabolic partial differential equations on $\R^N$ was used by Antocci and Prizzi \cite{Prizzi-tmna} and Prizzi \cite{Prizzi}. Recently a version of averaging principle
has been also obtained by the author in \cite{Cwiszewski-CEJM} where $A$ was a $C_0$ semigroup generator and $F$ a time periodic continuous perturbation.\\
\indent In this paper we look for a general averaging scheme in the abstract operator setting with time dependent $A$ and apply it to hyperbolic evolution equations. We shall prove a general principle, a version of which can be stated as follows (cf. Theorem \ref{25112009-1155} and Remark \ref{09072010-1159}).
\begin{Th}\label{09072010-1201}
Let $\{R^{(\lma)}(t,s)\}_{t\geq s\geq 0}$, $\lma>0$, be linear evolution systems on a separable Banach space $E$, corresponding to the problems
$$ \left\{
\begin{array}{l}
\dot u(t) = A(t/\lma)u(t), t>s,\\
u(s) = \bar u \in E.
\end{array} \right. $$
Suppose that\\
$(A1)$ \parbox[t]{138mm}{there are $M\geq 1$ and $\omega\in\R$ such that
$\|R^{(\lma)}(t,s)\| \leq Me^{\omega(t-s)}
\mbox{ if } t\geq s \geq 0$;}\\[1em]
$(A2)$ \parbox[t]{138mm}{there exists a $C_0$ semigroup $\{\widehat S(t)\}_{t\geq 0}$ of bounded linear operators on $E$
with the infinitesimal generator $\widehat A$ such that, for any $\bar u\in E$
and $t,s\geq 0$ with $t\geq 0$,
$$
\lim_{\lma\to 0^+, \, \bar v\to \bar u} R^{(\lma)} (t,s) \bar v  = \widehat S (t-s)\bar u
$$
uniformly with respect to $t$, $s$ from bounded intervals;}\\[1em]
$(A3)$ \parbox[t]{138mm}{$F:[0, +\infty) \times E \to E$ is Lipschitz on bounded subsets and has sublinear growth uniformly with respect to the second variable;}\\[1em]
$(A4)$  \parbox[t]{138mm}{ for each $\bar u\in E$, the set $\{ F(t, \bar u)\mid t\geq 0\}$ is relatively compact and
there is a locally Lipschitz mapping $\widehat F: E \to E$ such that, for any $\bar u\in E$ and $h>0$,
$$
\widehat F (\bar u) = \lim\limits_{T\to +\infty,\ \bar v\to\bar u}
\frac{1}{T}\int_{0}^{T} F(\tau+h,\bar v) \, \d\tau
$$
uniformly with respect to $h$.}\\[1em]
Then, for any $(\lma_n)$ in $(0,+\infty)$ and $(\bar u_n)$ in $E$ such that $\lma_n \to 0^+$ and $\bar u_n \to \bar u_0$ for some $\bar u_0\in E$,
the mild solutions $u_n:[0,+\infty)\to E$ of $(P_{\lma_n})$ satisfying $u_n (0)= \bar u_n$, $n\geq 1$, converge uniformly on bounded intervals
to the mild solution of the averaged problem
$$\left\{\begin{array}{l}
\dot u(t) = \widehat A u(t) + \widehat F (u(t)), \, t>0,\\
u(0) = \bar u_0.
\end{array}
\right.
$$
\end{Th}
\noindent The assumptions  $(A1)$ and $(A2)$ actually state that the averaging principle holds for the linear equation. Obviously, it is always the case if $A$ is independent of time and is an infinitesimal generator of a $C_0$ semigroup.
We will verify (A1) and (A2) for hyperbolic type linear evolution systems introduced by Kato -- see Theorem \ref{02032010-0959}.
Assumption (A3) and the separability of $E$  are to assure the existence of unique mild solutions for initial
value problems associated with $(P_\lma)$, the boundedness of solutions starting from bounded
sets and the relative compactness of semiorbits of relatively
compact sets (see $({\cal H}_1)-({\cal H}_3)$). Finally, $(A4)$ simply says that $F$ has the average $\widehat F$.
It is worth mentioning that $(A4)$ is fulfilled if $F$ is almost periodic with respect to time
(see \cite{Levitan-Zhikov}) and it is always the case when $F$ is time-periodic.
The obtained theorem generalizes those known in the literature -- see Remark
\ref{06072010-1905} adn besides the proof is rather straightforward.\\
\indent The paper is organized as follows. Section 2 is devoted to the general version of averaging principle
while in Section 3 we are concerned with its verification for abstract linear hyperbolic evolution systems.
Section 4 provides an example of application to first order hyperbolic partial differential equations.

\noindent {\bf Notation}\\
\indent By $\R$ we denote the field of real numbers; by $[x]$ we mean the integer (or floor) part of $x\in \R$.\\
\indent If $X$ is a metric space  and $B \subset X$, then
$\partial B$ and $cl B$ stand for the boundary of $B$ and the closure of $B$, respectively. If $x_0\in X$ and $r>0$, then
$B (x_0,r):=\{x\in M\mid d(x,x_0)<r  \}$.\\
\indent  If $E$ is a normed space, then by $\|\cdot\|$ we denote its norm.
If $V$ is another normed space then ${\cal L}(V,E)$ stands for the space of
all bounded linear operators with domain $V$ and values in $E$ with the operator norm denoted by $\| \cdot \|_{{\cal L}(V,E)}$ or simply $\| \cdot \|$
if no confusion may appear.

\section{General averaging principle}

Recall that a family of bounded linear operators $\{R(t,s)\!:\!E\!\to \!E\}_{t \geq s\geq 0}$
on a Banach space $E$ is an {\em evolution system} if and only
if $R(t,t)= I$, $R(t,s)R(s,r) = R(t,r)$, whenever $t\geq s\geq r\geq 0$, and
the mapping $(s,t)\mapsto R(t,s)\bar u$ is continuous for any $\bar u\in E$.
In this section we deal with general evolution systems, i.e. we do not indicate how they are generated. \\
\indent Evolution systems come up naturally in equations involving time-dependent families of linear operators.
Namely, if $\{ A (t) \}_{t\geq 0}$ is a family of linear operators in a Banach space $E$ satisfying suitable assumptions, then
for any $s\geq 0$ and $\bar u\in E$, the problem
$$\left\{
\begin{array}{l}
\dot u (t) = A(t)u(t), \quad t > s\\
u(s)=\bar u
\end{array}\right.
$$
admits a unique solution $u_{s,\bar u}:[s,+\infty)\to E$ (understood in an appropriate sense).
For instance, this is the case if $A(t)=A_0$, for each $t\geq 0$, with some $A_0$ being a generator of a $C_0$ semigroup of bounded linear operators on $E$,
as well as if $\{ A(t)\}_{t\geq 0}$ satisfies the so-called {\em parabolic} or {\em hyperbolic} conditions (see  e.g. \cite{Tanabe}, \cite{Pazy}
or \cite{Daners-book}). Moreover, the formula
$$
R(t,s)\bar u:= u_{s,\bar u}(t)  \mbox{ for } t\geq s
$$
defines an  evolution system $\{R(t,s)\}_{t\geq s \geq 0}$ on $E$.
In what follows we assume that cosidered evolution systems are generated by family time indexed families of operators in the above manner. We briefly
say that the evolution system $\{R(t,s)\}_{t\geq s \geq 0}$ is {\em determined by} or {\em correspond to} the family $\{A(t) \}_{t\geq 0}$.\\
\indent Let $\{ R(t,s) \}_{t \geq s \geq 0}$ be an arbitrary evolution system determined by a family $\{A(t):D(A(t))\to E\}_{t\geq 0}$
of linear operators in $E$. Consider the problem
$$
\left\{\begin{array}{l}
\dot u(t) = A(t)u(t)+F(t,u(t)),\, t\in [0,\omega),\\
u(0) = \bar u.
\end{array}
\right.
$$
where $\bar u\in E$, $\omega\in (0,+\infty]$ and $F:[0,+\infty)\times E\to E$ is a continuous mapping.
By a {\em mild solution} of the above problem we understand a continuous function $u:[0,\omega)\to E$  such that
$$
u(t)= R(t,0)\bar u + \int_{0}^{t}R(t,\tau)F(\tau,u(\tau))\,\d\tau \mbox{ for any } t\in [0,\omega).
$$
\indent We shall say that a family $\{R^{(\mu)}\}_{\mu\in P}$ of evolution systems, where $P$ is a metric space of parameters, is
{\em continuous} if, for any $\bar u\in E$ and $(\mu_n)$ in $P$ with $\mu_n\to \mu$,
$R^{(\mu_n)}(t,s)\bar u\to R^{(\mu)}(t,s) \bar u$ uniformly with respect to $t\geq s\geq 0$ from bounded intervals.

Now we pass to the averaging principle. For the sake of generality and future reference, we shall consider
its parameterized version. To this end we take families of operators $\{ A^{(\mu)}(t) \}_{t\geq 0}$, $\mu\in P$, where $P$ is a metric space of
parameters, determining corresponding evolution systems  $\{ R^{(\mu)}(t,s) \}_{t\geq s\geq 0}$, $\mu\in P$, on a Banach space $E$.
Assume that the family $\{R^{(\mu)}\}_{\mu\in P}$ is continuous and
that $F:[0,+\infty)\times E\times P \to E$ is a continuous mapping. Let families $\{
A^{(\mu,\lma)} (t)\}_{t\geq 0}$, $\mu\in P$, $\lma>0$, be defined by
$$
A^{(\mu,\lma)}(t):=A^{(\mu)}(t/\lma), \, t\geq 0
$$
and  $F^{(\mu,\lma)}:[0,+\infty)\times E\to E$, $\mu\in P$, $\lma>0$, by
$$
F^{(\mu,\lma)}(t,\bar u):=F(t/\lma, \bar u,\mu), \, t\geq 0, \,
\bar u\in E.
$$
The evolution system determined by $\{A^{(\mu, \lma)} (t)\}_{t\geq 0}$,  for $\mu\in P$, $\lma>0$, is denoted by
$\{R^{(\mu, \lma)} (t,s)\}_{t\geq s\geq 0}$.\\
\indent  We shall assume that the following conditions hold\\[1em]
$({\cal H}_1)$ \parbox[t]{140mm}{ for any $\bar u\in E$, $\mu\in P$ and $\lma>0$,  the problem
$$
\left\{
\begin{array}{l}  \dot u (t) =  A^{(\mu, \lambda)} u(t) +F^{(\mu,\lma)} (t, u(t)), \ t>0 \\
                              u (0) = \bar u \end{array} \right.
$$
admits a unique maximal mild solution $u(\cdot; \bar u, \mu, \lma): [0, \omega_{\bar u, \mu, \lambda})\to E$
with $\omega_{\bar u, \mu, \lambda} \in (0,+\infty]$;}\\[1em]
$({\cal H}_2)$ \parbox[t]{140mm}{given a bounded set $Q\subset E$, the sets
$F([0,+\infty) \times Q \times P)$ and $\{ u(t; \bar u, \mu, \lma) \mid t\in [0,\bar t],\, \bar u\in Q, \, \mu \in P, \, \lma>0 \}$, where $\bar t>0$ is such that $\bar t< \omega_{\bar u, \mu, \lma}$ for any $\bar u\in E$, $\mu\in P$ and $\lma>0$, are bounded;}\\[1em]
$({\cal H}_3)$ \parbox[t]{140mm}{if $Q_0\subset E$ and $P_0\subset P$
are relatively compact and $0<\bar t< \omega_{\bar u, \mu, \lma}$ for any $\bar u\in E$, $\mu\in P$ and $\lma>0$, then
$\{ u(\bar t; \bar u, \mu, \lma) \mid \bar u\in Q_0, \, \mu \in P_0, \, \lma>0 \}$
is relatively compact.}\\[1em]
$({\cal H}_4)$ \parbox[t]{140mm}{there are $M\geq 1$ and
$\omega\in\R$ such that, for any $\mu\in P$ and $\lma>0$,
$$
\|R^{(\mu,\lma)}(t,s)\| \leq M e^{\omega (t-s)} \mbox{ whenever } 0\leq s\leq t;
$$}\\[1em]
$({\cal H}_5)$ \parbox[t]{140mm}{ there are $C_0$ semigroups
$\widehat S^{(\mu)}=\{\widehat S^{(\mu)}(t):E\to E \}_{t\geq 0}$, $\mu\in P$, of bounded linear operators on
$E$ with the infinitesimal generators $\widehat A^{(\mu)}:D(\widehat A^{(\mu)})\to E$, $\mu\in P$, such
that, for any $t\geq 0$, $s\in [0,t]$, $\mu\in P$ and $\bar u\in E$,
$$
\lim\limits_{\lma \to 0^+,\ \bar v \to \bar u,\ \nu \to \mu} R^{(\nu, \lma)} (t,s)\bar v = \widehat S^{(\mu )}(t-s)\bar u
$$
and the convergence is uniform for $t$ and $s$ from bounded intervals;}\\[1em]
$({\cal H}_6)$ \parbox[t]{140mm}{$F: [0,+\infty)\times E\times P \to E$ is continuous uniformly with respect to the first variable, the set
$\{F(t,\bar u, \mu)\mid t\geq 0 \}$ is relatively compact for any $\bar u\in E$ and $\mu\in P$, and there is a continuous $\widehat
F: E\times P \to E$ such that, for any $\bar u\in E$, $\mu\in P$ and $h>0$,
\be\label{08042011-1443}
\widehat F (\bar u,\mu) = \lim\limits_{T\to +\infty,\ \bar v\to\bar u,\ \nu\to\mu}
\frac{1}{T}\int_{0}^{T} F(\tau+h,\bar v, \nu) \,\d\tau
\ee
where the convergence is uniform with respect to $h>0$;}\\[1em]
$({\cal H}_7)$ \parbox[t]{140mm}{for any  $\bar u\in E$ and $\mu \in P$ the averaged problem
$$\left\{\begin{array}{l}
\dot u (t)= \widehat A^{(\mu)} u(t) +\widehat F(u(t),\mu), \ t>0 \\
u(0)=\bar u
\end{array}\right.
$$
admits a unique maximal mild solution $\widehat u (\cdot; \bar u, \mu):[0,\widehat \omega_{\bar u,\mu})\to E$ with some $\widehat\omega_{\bar u, \mu}\in(0,+\infty]$.}
\begin{Rem}{\em $\mbox{ }$\\
\indent (i) Assumptions $({\cal H}_1)$ and $({\cal H}_7)$ are standard local existence properties, which hold if $F$ and $\widehat F$ are locally Lipschitz in the state variable.
We shall show in Proposition \ref{06032010-2035} that $({\cal H}_2)$
and $({\cal H}_3)$ hold for a large class of $F$.
Property $({\cal H}_4)$ is natural and is satisfied for example for the class of hyperbolic evolution systems considered in Section 3. \\
\indent (ii) Note that $({\cal H}_6)$ is a sort of an almost periodicity assumption (cf. \cite{Levitan-Zhikov}).
Moreover, it is always satisfied if $F$ is continuous and time periodic.\\
\indent (iii) Note that in $({\cal H}_5)$ we actually require that the averaging principle is true in the linear case. In Section 3 we shall prove it in the general hyperbolic case -- see Theorem \ref{02032010-0959}.
}\end{Rem}

\begin{Th} {\em (Abstract averaging principle)} \label{25112009-1155}
Let $({\cal H}_1)-({\cal H}_7)$ be satisfied.
If $(\bar u_n)$ in  $E$, $(t_n)$ in $[0,+\infty)$, $(\mu_n)$ in $P$
and $(\lma_n)$ in $(0,+\infty)$ are such that $\bar u_n \to \bar u_0$, $t_n\to t_0$, $\mu_n\to \mu_0$, $\lma_n\to 0^+$
as $n\to +\infty$, for some $t_0\geq 0$, $\bar u_0\in E$ and $\mu_0\in P$, and $t_n \leq \bar t < \omega_{\bar u_n, \mu_n, \lma_n}$ for some $\bar t>0$ and each
$n\geq 1$, then
$$
u(t_n; \bar u_n,\mu_n, \lma_n)  \to \widehat u (t_0; \bar u_0,\mu_0) \mbox{ in } E \mbox{ as } n\to\infty.
$$
Moreover
\be\label{27042010-1105}
\ \ \max\{
\|u(t; \bar u_n,\mu_n, \lma_n)\!-\!\widehat u (t; \bar  u_0,\mu_0)\|\  |\ t\in [0, \bar t]\} \longrightarrow 0
\mbox{ as  } n\!\to\! +\infty.
\ee
\end{Th}
To prove it we shall need three auxillary facts.
\begin{Lem}\label{05122009-1806} {\em (See \cite[Proposition 3.1]{Cwiszewski-Kokocki-1})}
Suppose that $\{R_n(t,s)\}_{t\geq s\geq 0}$, $n\geq 1$, are evolution systems on $E$ with $M\geq 0$ and $\omega\in\R$ such that
$$
\|R_n (t,s)\| \leq Me^{\omega(t-s)}, \mbox{ for any } t,s\geq 0 \mbox{ with } t\geq s,
$$
and there is an evolution system $\{R(t,s):E\to E\}_{t\geq s\geq 0}$ such that
$$
\lim_{n\to +\infty }R_n (t,s) \bar u = R(t,s)\bar u \mbox{ for all }\bar u\in E.
$$
Let $\{\bar u_n\}_{n\geq 1 }\subset E$ be relatively compact
and $\{w_n\}_{n\geq 1} \subset L^1([0,l],E)$ be uniformly integrable {\em(\footnote{i.e., for any $\eps>0$, there is $\delta>0$ such that,
for any measurable $J\subset [0,l]$ with the Lebesgue measure $\mu(J) \leq\delta$ and any $n\geq 1$, $\int_{J} \|w_n(t)\| dt\leq \eps$})}.
Put $u_n:[0,l]\to E$, $n\geq 1$, by
$$
u_n(t):=R_n(t,0)\bar u_n + \int_{0}^{t} R_n (t,s)w_n(s) \, \d s, \ t\in [0,l].
$$
Then the following conditions are equivalent\\
\indent {\em (a)} $\{ u_n(t) \}_{n\geq 1}$ is relatively compact for a.e. $t\in [0,l]$;\\
\indent {\em (b)}  \parbox[t]{140mm}{$\{ u_n \}_{n\geq 1}$ is relatively compact in the space $C([0,l],E)$ (with the uniform convergence norm).}
\end{Lem}
\begin{Lem} \label{31052010-1559}
Let $({\cal H}_6)$ be satisfied and $Q\subset E$ be compact.
Then, for any
$(T_n)$ in $(0,+\infty)$ with $T_n\to +\infty$ and $(\mu_n)$ in $P$ with $\mu_n\to \mu_0$,
the convergence
$$
\lim_{n\to +\infty} \frac{1}{T_n}\int_{0}^{T_n} F(\tau+h, \bar w,\mu_n)\, \d\tau = \widehat F (\bar w, \mu_0 )
$$
is uniform with respect to $\bar w\in Q$ and $h>0$.
\end{Lem}
{\bf Proof:} Put $P_0 := cl \, \{\mu_n \mid n \geq 1\}$. Take
an arbitrary $\eps>0$. By the compactness of $Q$ and $P_0$ there
is a set $\{ \bar w_k\}_{k=1}^{n_\eps}\subset Q$ with
$\delta_{1},\ldots, \delta_{n_\eps} \in (0,\eps)$ such that
$$
Q \subset \bigcup_{k=1}^{n_\eps}  B ( \bar w_k , \delta_{n_k} )
$$
and, for any $\tau>0$, $\mu \in P_0$ and $\bar w \in B(\bar w_k,\delta_{k})$, $k\in \{1,\ldots, n_\eps\}$,
\begin{eqnarray*}
& \| \widehat F (\bar w, \mu) - \widehat F(\bar w_k,\mu)\| < \eps/3 &\\  
& \| F(\tau, \bar w, \mu) - F(\tau, \bar w_k, \mu) \|<\eps/3. & 
\end{eqnarray*}
Due to (\ref{08042011-1443}), there exists $n_0\geq 1$ such that, for any $n\geq n_0$,  $k=1,\ldots, n_\eps$ and $h>0$,
$$\left\| \frac{1}{T_n} \int_{0}^{T_n} F(\tau + h, \bar w_k, \mu_n)\, \d\tau - \widehat F (\bar w_k, \mu_0) \right\| <\eps/3.$$
Taking $n\geq n_0$, $\bar w\in Q$ and $h>0$, we get
$\bar w\in B(\bar w_k, \delta_k)$ for some $k=1, \ldots, n_\eps$ and,
consequently,
\begin{eqnarray*}
& & \left\| \frac{1}{T_n} \int_{0}^{T_n} F(\tau+h, \bar w, \mu_n)\, \d\tau-\widehat F(\bar w,\mu_0) \right\| \\
& & \ \ \ \ \   \leq \frac{1}{T_n} \int_{0}^{T_n} \left\|  F(\tau+h, \bar w, \mu_n)-F(\tau+h,\bar w_k, \mu_n) \right\| \, \d\tau\\
& & \ \ \ \ \ \ \ \ \ \ \ \ \  +\left\| \frac{1}{T_n} \int_{0}^{T_n}F(\tau+h, \bar w_k, \mu_n)\, \d\tau  - \widehat F(\bar w_k, \mu_0)  \right\|+
\left\|\widehat F(\bar w_k,\mu_0)-\widehat F(\bar w, \bar \mu_0) \right\| \\
& & \leq \eps/3 + \eps/3 + \eps/3 = \eps,
\end{eqnarray*}
which  completes the proof. \hfill $\square$
\begin{Lem}\label{01032010-2300}
Assume that $({\cal H}_4)$ and $({\cal H}_5)$ hold. Let $(T_n)$ be a sequence in $(0,+\infty)$, $(k_n)$ a sequence of positive integers,
$(\lma_n)$ a sequence in $(0,+\infty)$ and $(\mu_n)$ in $P$ such that
$T_n \to +\infty$, $k_n\to \infty$, $\lma_n \to 0$, $k_n \lma_n T_n \to t$ for some $t>0$, $\mu_n \to \mu_0$ for some
$\mu_0\in P$, as $n\to +\infty$,  and $k_n \lma_n T_n \leq t$ for almost all integers $n\geq 1$.
Then, for any continuous function $w:[0,t] \to E$,
$$
\lma_n T_n \sum_{k=0}^{k_n-1} R^{(\mu_n,\lma_n)}(k_n \lma_n T_n, k\lma_n T_n) w(k\lma_n T_n) \to \int_{0}^{t} \widehat S^{(\mu_0)}(t-s) w(s) d s  \ \mbox{ as } n\to +\infty.
$$
\end{Lem}
{\bf Proof:}  First note that, for any $n\geq 1$,
$$
\lma_n T_n \sum_{k=0}^{k_n-1} R^{(\mu_n,\lma_n)}(k_n \lma_n T_n, k\lma_n T_n) w(k\lma_n T_n)   = \sigma_{1,n} +
\sigma_{2,n} + \sigma_{3,n}
$$
where
\begin{eqnarray*}
\sigma_{1,n}& := & \lma_n T_n \sum_{k=0}^{k_n-1}  [R^{(\mu_n,\lma_n)}(k_n \lma_n T_n, k\lma_n T_n) - \widehat S^{(\mu_0)}( k_n \lma_n T_n - k\lma_n T_n)]w(k\lma_n T_n),\\
\sigma_{2,n}& := & \lma_n T_n \sum_{k=0}^{k_n-1}  [\widehat S^{(\mu_0)}( k_n \lma_n T_n - k\lma_n T_n)- \widehat S^{(\mu_0)}( t- k\lma_n T_n)]w(k\lma_n T_n), \\
\sigma_{3,n}& := & \lma_n T_n \sum_{k=0}^{k_n-1}  \widehat S^{(\mu_0)}( t- k\lma_n T_n)w(k\lma_n T_n).
\end{eqnarray*}
It follows from $({\cal H}_5)$ and  the compactness of $w([0,t])$ that, for any $\eps>0$,
there exists $n_0\geq 1$ such that, for all $n\geq n_0$,
$\bar w\in w([0,t])$ and  $s',s\in [0,t]$ with  $s'\geq s$,
$$
\| R^{(\mu_n,\lma_n)} (s',s) \bar w - \widehat S^{(\mu_0)}(s'-s) \bar w\| \leq \eps.
$$
This implies $\sigma_{1,n}\to 0$ as $n\to +\infty$, since $k_n\lma_n T_n \to t$.
Again, by the compactness of $w([0,t])$ and the strong continuity of the semigroup $\widehat S^{(\mu_0)}$, we gather that
$\sigma_{2,n}\to 0$ as $n\to +\infty$. Finally, by the continuity of $[0,t]\ni s\mapsto \widehat S^{(\mu_0)} (t-s)w(s)$ and the fact that
$\lma_n T_n \to 0$ as $n \to +\infty$, we infer that
$$
\sigma_{3,n}\to \int_{0}^{t} \widehat S^{(\mu_0)} (t-s) w(s) \, \d s,
$$
which ends the proof. \hfill $\square$

\noindent {\bf Proof of Theorem \ref{25112009-1155}:}
Let  $R_n:=R^{(\mu_n, \lma_n)}$ and $u_n:[0,\bar t]\to E$, $n\geq 1$, be given by
\be\label{04012010-1825}
u_n(t):=u(t; \bar u_n, \mu_n, \lma_n) \mbox{ for } t\in [0,\bar t],  \ n\geq 1. \ee
Obviously
\be\label{1927-19062010}
u_n(t)= R_n (t,0) \bar u_n + \int_{0}^{t} R_n (t,s)w_n(s) \, \d s \mbox{ for  } t\in [0,\bar t]
\ee
with $w_n:[0,\bar t]\to E$ defined as $w_n (s):= F^{(\mu_n, \lma_n)} (s, u(s;\bar u_n, \mu_n,\lma_n))$ for $s\in [0,\bar t]$.
By $({\cal H}_2)$, there are $C_1, C_2 >0$
such that
\be\label{08122009-0100} \| u_n(t) \|  \leq C_1 \ \ \
\mbox{ and } \ \ \ \|w_n(t)\| \leq C_2 \mbox{ for all } t\in
[0,\bar t] \mbox{ and } n\geq 1.
\ee
In the rest of the proof we shall argue as follows: we take any subsequence of $(u_n)$, denote it again by $(u_n)$ and show that it contains a subsequence
converging to $\widehat u (\cdot; \bar u_0,\mu_0)_{|[0,\bar t]}$ in the space $C([0,\bar t],E)$;
having this we will conclude that the original $(u_n)$ converges
uniformly to $\widehat u (\cdot; \bar u_0,\mu_0)$ on $[0,\bar t]$ and the assertion will follow.\\
\indent Start with an observation that, in view of $({\cal H}_3)$,
for each $t \in [0,\bar t]$, the set $\{u_n(t) \}_{n\geq 1}$ is relatively compact.  And due to $({\cal H})_4$, $({\cal H}_5)$ and  Lemma \ref{05122009-1806}, $(u_n)$ contains a subsequence converging uniformly on $[0, \bar t]$ to some $\widehat v$. Denote that subsequence again by $(u_n)$.
It follows directly from $({\cal H}_5)$ that
\be\label{01032010-1304}
R_n (t,0) \bar u_n  = R^{(\mu_n,\lma_n)} (t,0) \bar u_n \to \widehat S^{(\mu_0)} (t) \bar u_0 \mbox{ as } n\to +\infty
\ee
uniformly with respect to $t\in [0,\bar t]$.\\
\indent Now take an arbitrary  $t\in (0,\bar t]$ and any sequences $(T_n)$ in $(0,+\infty)$
and $(k_n)$ of positive integers such that $T_n \to +\infty$, $k_n \to +\infty$, $k_n \lma_n T_n \to t$, as $n\to +\infty$,
and $k_n \lma_n T_n \leq t$ for any $n\geq 1$ (e.g. $T_n:=\lma_n^{-1/2}$, $k_n:= [t/\lma_n T_n]$). First observe that
$$
\int_{0}^{t} R_n (t,s)w_n(s)\,  \d s = I_{1,n} + I_{2,n} + I_{3,n}
$$
where
\begin{eqnarray*}
I_{1,n} &:=& \int_{0}^{k_n\lma_n T_n} R_n (k_n\lma_n T_n,s)w_n(s)\, \d s,\\
I_{2,n} &:=& \int_{0}^{k_n\lma_n T_n}
 \ \ \ \  (R_n(t,s) - R_n(k_n\lma_n T_n,s)) w_n (s)\, \d s, \\
I_{3,n} &:=& \int_{k_n\lma_n T_n}^{t} \!\!\!\! R_n (t,s)w_n(s)\, \d s. \\
\end{eqnarray*}
It is immediate to see that, by $({\cal H}_4)$ and (\ref{08122009-0100}), $I_{3,n} \to 0$ $\mbox{as } n\to +\infty$.\\
\indent To deal with $(I_{2,n})$, we claim that
\be\label{01042011-1441}
\mbox{ the set $\widetilde Q:=\{w_n (s)\mid s\in [0,t], n\geq 1 \}$ is relatively compact}.
\ee
Indeed, to see this take any sequence $(\bar w_k)$ in $\widetilde Q$. Then, for each $k\geq 1$, there are an integer $n_k\geq 1$ and $s_k \in [0,t]$
such that
$$
\bar w_k = w_{n_k}(s_k) = F(s_k / \lma_{n_k}, u_{n_k}(s_k), \mu_{n_k}).
$$
We may assume that $s_k \to s$ as $k\to +\infty$ and, by the uniform convergence of $(u_n)$, that $u_{n_k} (s_k)\to \widetilde u$
for some $\widetilde u\in E$. In view of $({\cal H}_6)$, for any $\eps>0$ one finds $n_0\geq 1$
such that, for all $n\geq n_0$,
$$
\|\bar w_k - F(s_k/\lma_{n_k}, \widetilde u, \mu_0) \|=\| F(s_k/\lma_{n_k}, u_{n_k}(s_k), \mu_{n_k}) - F(s_k/\lma_{n_k}, \widetilde u, \mu_0) \| < \eps.
$$
Now, since, by $({\cal H}_6)$, the set $\{F(s, \widetilde u, \mu_0) \mid t\geq 0\} $ is relatively compact, it follows that $(\bar w_k)$ contains a convergent subsequence, which proves (\ref{01042011-1441}).\\
\indent It can be easily seen that $({\cal H}_5)$ along with (\ref{01042011-1441}) implies that,
for any $s',s \in [0,t]$ with $s'\geq s$ and $\bar w \in \widetilde Q$,
\be\label{1952-19062010} R_n (s',s) \bar w \to \widehat
S^{(\mu_0)} (s'-s) \bar w \mbox{ as  } n\to +\infty \mbox{ uniformly with respet to } s,s' \mbox{ and } \bar w.
\ee
In consequence, for any $\eps>0$ there exists $n_1\geq 1$ such that, for all $n\geq n_1$,
\begin{eqnarray}\label{01042011-1500}
\| R_n(t,s) w_n (s) - \widehat S^{(\mu_0)} (t-s) w_n(s)\|  & \leq & \eps/3t, \mbox{ for } s\in [0,t], \\
 \|\widehat S^{(\mu_0)} (k_n \lma_n T_n -s) w_n(s)- R_n (k_n\lma_n T_n,s) w_n (s) \|& \leq & \eps/3t, \mbox{ for } s\in [0, k_n \lma_n T_n]. \label{01042011-1501}
\end{eqnarray}
Moreover, by the strong conitnuity of the semigroup $\widehat  S^{(\mu_0)}$ and the relative compactnes of $\widetilde Q$,
there is $n_0\geq n_1$ such that, for all $n\geq n_0$ and $s\in [0,t]$,
\be\label{01042011-1502}
\| \widehat  S^{(\mu_0)} (t-s) w_n(s) - \widehat S^{(\mu_0)} (k_n \lma_n T_n -s) w_n(s)\| <\eps/3t.
\ee
Since it follows from (\ref{01042011-1500}), (\ref{01042011-1501}) and (\ref{01042011-1502}) that, for all $n\geq n_0$ and $s\in [0,t]$,
$$
\|(R_n(t,s) - R_n(k_n\lma_n T_n,s)) w_n (s)\| \leq \eps/3t + \eps/3t + \eps/3t = \eps/t
$$
we infer, for $n\geq n_0$, $\| I_{2,n} \| \leq k_n \lma_n T_n (\eps/t) \leq \eps$,
i.e., $I_{2,n} \to 0$ $\mbox{as } n\to +\infty$.\\
\indent Now our aim is to show that
$$
I_{1,n} \to \int_{0}^{t} \widehat S^{(\mu_0)}(t-s) \widehat F (\widehat v (s),\mu_0) d s \mbox{ as } n\to +\infty.
$$
To this end observe that, for sufficiently large $n\geq 1$,
\begin{eqnarray*}
I_{1,n} &=& \int_{0}^{k_n\lma_n T_n} R_n (k_n\lma_n T_n,s)w_n(s)\, \d s
= \sum_{k=0}^{k_n-1} \int_{k\lma_n T_n}^{(k+1)\lma_n T_n}
\!\!\!\! R_n(k_n\lma_n T_n,s)w_n(s) \, \d s \\
&= &\sum_{k=0}^{k_n-1} \int_{0}^{T_n}
\!\!\!\! R_n(k_n\lma_n T_n,k\lma_n T_n + \lma_n\tau )w_n(k\lma_n T_n +\lma_n\tau )\lma_n \, \d\tau \\
& = & J_{1,n} + J_{2,n} + J_{3,n} + J_{4,n}\\
\end{eqnarray*}
with
\begin{eqnarray*}
\!\!\!\! J_{1,n}&\!\!\!\! :=& \!\!\!\!\!\!\sum_{k=0}^{k_n-1} \!\int_{0}^{T_n} \!\! \left[R_n(k_n\lma_n T_n, \! k\lma_n T_n\!
+\!\lma_n \tau)-R_n(k_n\lma_n T_n, \! k\lma_n T_n)\right] w_n (k\lma_n T_n\!+\!\lma_n\tau)\lma_n\, \d\tau\\
\!\!\!\! J_{2,n}&\!\!\!\! :=& \!\!\!\!\!\!\sum_{k=0}^{k_n-1} \! R_n(k_n\lma_n T_n,  k\lma_n T_n)
\int_{0}^{T_n} \!\left(w_n(k\lma_n T_n+\lma_n \tau)-
   F(k T_n+\tau,\widehat v(k\lma_n T_n),\mu_n)\right)\lma_n\, \d\tau\\
\!\!\!\! J_{3,n} &\!\!\!\!:= & \!\!\!\!\!\! \sum_{k=0}^{k_n-1}\!\! R_n(k_n\lma_n T_n,  k\lma_n T_n)
\frac{1}{T_n}\int_{0}^{T_n} \!\!\!
(F(k T_n\!+\!\tau,\widehat v (k\lma_n T_n),\mu_n)-
     \widehat F (\widehat v (k\lma_n T_n),\mu_0) ) \lma_n T_n\, \d\tau\\
\!\!\!\! J_{4,n} &\!\!\!\! := & \!\!\!\!\!\! \sum_{k=0}^{k_n-1} \!\! R_n(k_n\lma_n T_n,  k\lma_n T_n)
\widehat F (\widehat v (k\lma_n T_n),\mu_0) )\lma_n T_n.
\end{eqnarray*}
First we note that, for all $n\geq 1$, $\tau \in [0,T_n]$ and $k=0,1, \ldots, k_n-1$,
\begin{eqnarray*}
& & \|[R_n(k_n\lma_n T_n, \! k\lma_n T_n\! +\!\lma_n \tau)\!-\!R_n(k_n\lma_n T_n, \! k\lma_n T_n) ] w_n (k\lma_n T_n\!+\!\lma_n\tau)\|\\
& & \ \ \ \ \leq \| [R_n(k_n\lma_n T_n, \! k\lma_n T_n\! +\!\lma_n \tau)\!-\!\widehat S^{(\mu_0)}
(k_n\lma_n T_n \!-\! k\lma_n T_n\! -\!\lambda_n \tau)] w_n (k\lma_n T_n\!+\!\lma_n\tau) \|\\
& & \ \ \ \ \ \ \  + \| [\widehat S^{(\mu_0)} (k_n\lma_n T_n \!-\! k\lma_n T_n\! -\!\lambda_n \tau) - \widehat S^{(\mu_0)} (k_n\lma_n T_n \!-\! k\lma_n T_n)]
w_n (k\lma_n T_n\!+\!\lma_n\tau)\|\\
& & \ \ \ \ \ \ \  + \|[\widehat S^{(\mu_0)} (k_n\lma_n T_n \!-\! k\lma_n T_n) -  R_n(k_n\lma_n T_n, \! k\lma_n T_n)  ] w_n (k\lma_n T_n\!+\!\lma_n\tau)\|.
\end{eqnarray*}
Hence, in view of the uniform convergence in (\ref{1952-19062010}) and
the uniform equicontinuity of functions $[0,t] \ni s\mapsto \widehat S^{(\mu_0)} (s)\bar w$, $\bar w\in \widetilde Q$, we deduce that, for any $\eps>0$, there is $n_0\geq 1$, such that for  all  $n\geq n_0$, $\tau \in [0,T_n]$ and $k=0,1, \ldots, k_n-1$,
$$
\|[R_n(k_n\lma_n T_n, \! k\lma_n T_n\! +\!\lma_n \tau)\!-\!R_n(k_n\lma_n T_n, \! k\lma_n T_n) ] w_n (k\lma_n T_n\!+\!\lma_n\tau)\| < \eps/t
$$
and
$\|J_{1,n}\| < k_n \lma_n T_n (\eps/t) \leq \eps.$
This means that $J_{1,n}\to 0$ as $n\to + \infty$.\\
\indent As for $J_{2,n}$, take any $\eps>0$ and note that, by the uniform convergence of $(u_n)$, the set
$$
\widetilde Q':=cl \left(\bigcup_{n\geq 1} u_n([0,t])\right)
$$
is compact and $\widehat v ([0,t])\subset \widetilde Q'$. Therefore, by $({\cal H}_6)$, we find $\eta >0$ such that,
for any $\tau\geq 0$, $\mu \in cl \{\mu_n\,|\, n\geq 1 \}$, $\bar v_1, \bar v_2\in \widetilde Q'$
with  $\|\bar v_1-\bar v_2 \|\leq \eta$,
$$
\|F(\tau, \bar v_1, \mu) - F(\tau,\bar v_2, \mu)\| \leq \eps/ t Me^{|\omega| t}.
$$
Furthermore, again by the uniform convergence of $(u_n)$, there is $\delta>0$ and $n_1\geq 1$ such that, for any $n\geq n_1$,
$$
\|u_n(\tau_1)-\widehat v(\tau_2)\|<\eta \, \mbox{ for any } \, \tau_1, \tau_2\in [0,t] \mbox{ with } |\tau_1-\tau_2|<\delta.
$$
Let $n_0\geq n_1$ be such that $\lma_n T_n <\delta$ for each $n\geq n_0$. Then, for any $n\geq n_0$, $\tau\in [0,T_n]$
and $k\in \{0,\ldots, k_n-1\}$,
$$
\| u_n (k\lma_n T_n + \lma_n \tau ) -\widehat  v (k \lma_n T_n)\| \leq \eta,
$$
and, consequently, for any $n\geq n_0$, $\tau\in [0,T_n]$ and $k=0,1,\ldots, k_n-1$,
\begin{eqnarray*}
& & \|w_n(k\lma_n T_n\!+\!\lma_n\tau) \!-\! F(k T_n \!+\!\tau,\widehat v(k\lma_n T_n),\mu_n)\|\\
& & \ \ \ \ \ \ \ \ =  \|F(k T_n \!+\!\tau, u_n ( k\lma_n T_n \!+\!\lma_n \tau ),\mu_n )\!-\! F(k T_n \!+\!\tau, \widehat v (k\lma_n T_n),\mu_n)\|
\leq \eps/t Me^{|\omega|t}.
\end{eqnarray*}
Therefore, for any $n\geq n_0$, one has
\begin{eqnarray*}
\|J_{2,n}\| \!\!\!\! &\leq\!\!\!\!& Me^{|\omega|t}\lma_n \sum_{k=0}^{k_n-1} \int_{0}^{T_n} \!\!\!
 \| w_n ( k\lma_n T_n\!+\! \lma_n\tau)\!-\!
F(k T_n + \tau, \widehat v (k\lma_n T_n),\mu_0)\| \, \d\tau \\
\!\!\!\!&<\!\!\!\!& Me^{|\omega|t}  \lma_n k_n T_n (\eps/ t M e^{|\omega|t})= k_n \lma_n T_n \eps/ t \leq \eps,
\end{eqnarray*}
which yields $J_{2,n}\to 0$ as $n\to +\infty$.\\
\indent  Next observe that, in view of Lemma \ref{31052010-1559},
$$
\lim_{n\to + \infty} \frac{1}{T_n}\int_{0}^{T_n} F(h+ \tau,\bar w,\mu_n)\, \d\tau = \widehat F (\bar w,\mu_0)
$$
uniformly with respect to $\bar w \in \widetilde Q'$ and $h>0$.
Hence
$$
\|J_{3,n}\|\leq Me^{|\omega|t} \lma_n \sum_{k=0}^{k_n-1} \left\|\frac{1}{T_n}\int_{0}^{T_n} \!\! F(k T_n\!+\!\tau,\widehat v(k\lma_n T_n),\mu_n)\d\tau
-     \widehat F (\widehat v (k\lma_n T_n),\mu_0) \right\|
$$
and we find that $J_{3,n}\to 0$ as $n\to +\infty$.\\
\indent Finally, due to Lemma \ref{01032010-2300},
$$
J_{4,n} \to \int_{0}^{t} \widehat S^{(\mu_0)} (t-s) \widehat F (\widehat v(s), \mu_0)\,  \d s.
$$
Summing up, we have already showed that
$$
\int_{0}^{t} R_n (t,s)w_n (s)\, \d s \to \int_{0}^{t} \widehat S^{(\mu)} (t-s) \widehat F (\widehat v(s), \mu_0)\, \d s
$$
(after passing to a subsequence). Thus, letting $n\to +\infty$ in
(\ref{1927-19062010}), we arrive at
$$
\widehat v(t)=\widehat S^{(\mu_0)}(t)\widehat v(0) + \int_{0}^{t} \widehat S^{(\mu_0)}(t-s) \widehat F(\widehat v(s),\mu_0)\, \d s, \mbox{ for any }
t\in [0,\bar t].
$$
In particular, $\widehat v$ is a mild solution of $\dot u\! = \!\widehat A^{(\mu_0)} u \!+\!\widehat F(u,\mu_0)$ and, in view of $({\cal H}_7)$,
$\widehat v = \widehat u (\cdot; \bar u_0, \mu_0)_{|[0,\bar t]}$. Hence the original sequence $(u_n)$ converges uniformly to $\widehat u (\cdot; \bar u_0, \mu_0)$ on $[0,\bar t]$, which together with (\ref{04012010-1825}), implies $u_n(t_n;\bar u_n, \mu_n, \lma_n)\!\to\! \widehat u (t_0; \bar u_0,\mu_0)$ and as well as (\ref{27042010-1105}).\hfill$\square$
\begin{Rem} \label{06072010-1905} {\em
(a) For single-valued $F$, Theorem \ref{25112009-1155} is an
extension of \cite[Theorem 5.4.1]{Ka-Obu-Zec} in a few aspects:
firstly, the periodicity assumption of $F$ is relaxed; secondly,
the case where $A$ depends on time is included; thirdly, the
sublinear growth condition is not required and, finally, no
explicit compactness assumptions on $F$ are imposed (which may be of importance, see e.g. \cite{Cwiszewski-CEJM}
where one could not apply \cite[Theorem 5.4.1]{Ka-Obu-Zec}).\\
\indent (b) Note also that Theorem 2.2 covers \cite[Theorem
2.4]{Cwiszewski-CEJM} as well.}\end{Rem}
We end this section with
an example showing the availability of assumptions $({\cal H}_2)$
and $({\cal H}_3)$ and derive Theorem \ref{09072010-1201}.
\begin{Prop}\label{06032010-2035}
Let $E$ be a separable Banach space and $\{A ^{(\mu)} (t)\}_{t\geq 0 }$, $\mu\in P$, be such that
the family $\{ R^{(\mu)}\}_{\mu\in P}$ of the corresponding evolution systems is continuous and satisfies $({\cal H}_4)$
and $({\cal H}_5)$.
Assume that a continuous $F:[0,+\infty) \times E\times P\to E$ satisfies the following conditions
\be\label{09072010-1216}
\mbox{\parbox[t]{130mm}{for any $\bar v\in E$, there exist $L>0$ and
$\delta >0$ such that
$$
\|F(t,\bar v_1,\mu)\!-\!F(t,\bar v_2,\mu)\|\leq L\|\bar v_1\! -\! \bar v_2\|\ \  \mbox{ for any  $t\geq 0$, $\bar v_1,\bar v_2\in B(\bar v,\delta)$, $\mu\in P;$}
$$}}
\ee
\be\label{1916-19062010}
\mbox{\parbox[t]{130mm}{there is $c>0$ such that
$$
\| F (t,\bar v,\mu ) \| \leq  c (1+\|\bar v\|)   \ \ \mbox{ for any  } (t,\bar v, \mu)\in [0,+\infty)\times E\times P;
$$}}
\ee
and
\be\label{06032010-2111}
\mbox{\parbox[t]{130mm}{ for any $C>0$ there is $k \geq 0$ such that
$$
\beta (F([0,+\infty)\!\times\! Q \!\times\! P )) \leq k \cdot \beta (Q) \mbox{ for any  } Q \subset B(0,C)
$$
where $\beta$ stands for the Hausdorff or Kuratowski measure of noncompactness {\em(}see e.g. {\em \cite{Deimling}} or {\em \cite{Ka-Obu-Zec}).}}}
\ee
\indent Then  $({\cal H}_1)-({\cal H}_3)$ are satisfied.
\end{Prop}
We shall use the following general properties involving the measures of noncompactness.
\begin{Lem}\label{13122009-2039} {\em (see \cite[Prop. 9.3]{Deimling} or \cite{Ka-Obu-Zec})}
Let $E$ be a separable Banach space, $W\subset L^1([0,l],E)$,
$l>0$, be countable and integrably bounded {\em(}i.e. there exists
$c\in L^1([0,l])$ such that $\|w(t)\|\leq c(t)$  for all $w\in W$
and a.e. $t\in [0,l]${\em)} and $\phi:[0,l]\to\R$ be given by
$\phi (t):=\beta(\{u(t)|\, u\in W \})$. Then $\phi\in L^1 ([0,l])$
and
$$
\beta \left(\left\{ \int_{0}^{l} u(\tau) \d \tau\,|\, u \in
W\right\} \right) \leq \int_{0}^{l} \phi(\tau) \d \tau.
$$
\end{Lem}
\begin{Lem}\label{13122009-2040} {\em (see \cite[Lemma 5.4]{Cwiszewski-Kokocki})}
Let $T_n:E\to E$, $n\geq 1$, be bounded linear operators on a Banach space $E$ such that, for any $\bar u \in E$, $(T_n \bar u)$ is a
Cauchy sequence. Then, for any bounded set $\{\bar u_n\}_{n\geq 1}
\subset E$,
$$\beta\left( \{T_n \bar u_n\}_{n\geq 1} \right) \leq \left( \limsup_{n\to +\infty}
\|T_n\|\right) \beta\left( \{\bar u_n\}_{n\geq 1}\right).$$
\end{Lem}
\noindent {\bf Proof of Proposition \ref{06032010-2035}:} It is standard to see that the local Lipschitzianity of $F$ implies
the local existence, i.e., $({\cal H}_1)$ holds. It can be also easily deduced that  the sublinear growth and the Gronwall inequality yield $({\cal H}_2)$.\\
\indent In order to verify $({\cal H}_3)$, take any $\{ \lma_n\}_{n\geq 1} \subset (0,+\infty)$ and relatively compact sets
$\{ \bar u_n \}_{n\geq 1}\subset E$, $\{ \mu_n \}_{n\geq 1} \subset P$ and suppose that $0 < t < \omega_{\bar u_n, \mu_n, \lma_n}$ for any $n\geq 1$.
Define $\phi: [0,t] \to \R$ by
$$
\phi(\tau):= \beta ( \{ u (\tau; \bar u_n, \mu_n, \lma_n ) \}_{n\geq 1} )
$$
and $w_n:[0,t]\to E$ by $w_n (\tau):= F(\tau/\lma_n, u (\tau;\bar u_n, \mu_n, \lma_n),\mu_n)$, $\tau\in [0,t]$.
By $({\cal H}_2)$, there is $C>0$ such that $\sup_{n\geq 1}\sup_{\tau \in [0,t]}\|w_n(\tau)\| < C$,
and, in view of Lemma \ref{13122009-2039}, $\phi$ is integrable.
Taking into account $({\cal H}_4)$ and $({\cal H}_5)$ and using Lemmas \ref{13122009-2040} and \ref{13122009-2039}, for any $r\in [0,t]$, one gets
\begin{eqnarray*}
\phi(r)& \leq&  \beta(\{R^{(\mu_n,\lma_n)}(r,0)\bar u_n\}_{n\geq 1}) + \beta\left(\left\{\int_{0}^{r} R^{(\mu_n, \lma_n)} (r,\tau) w_n(\tau)\,\d\tau\right\}_{n\geq 1}     \right)\\
& \leq & Me^{|\omega|t} \beta(\{\bar u_n \}_{n\geq 1}) + \int_{0}^{r} \beta (\{R^{(\mu_n, \lma_n)}(r,\tau)w_n(\tau)\}_{n\geq 1})\, \d\tau\\
& \leq &  Me^{|\omega|t} \int_{0}^{r} \beta (\{w_n(\tau)\}_{n\geq 1})\, \d\tau.
\end{eqnarray*}
Hence, by use of (\ref{06032010-2111}), there is $k>0$ such that
$$
\phi(r)\leq  Me^{|\omega|t} k \int_{0}^{r} \phi(\tau) \, \d\tau,
$$
which implies $\phi(r)=0$ and completes the proof of $({\cal H}_3)$. \hfill $\square$

\begin{Rem} \label{09072010-1159} {\em
Note that Theorem \ref{09072010-1201} is a consequence of Theorem \ref{25112009-1155}.
Indeed, the Lipshitzianity on bounded sets implies that (\ref{06032010-2111}) holds and, if $\{A(t)\}_{t\geq 0}$ satisfy $(A1)$ and $(A2)$,
then, in view of Proposition \ref{06032010-2035}, assumptions $({\cal H}_1)-({\cal H}_3)$ are fulfilled.}
\end{Rem}

\section{Averaging for linear hyperbolic evolution systems}

Let $V$ be a Banach space densely and continuously embedded into a Banach space $E$.
If a linear operator $A:D(A)\to E$ generates a $C_0$ semigroup $\{S_A(t)\}_{t\geq 0}$ of bounded linear operators on $E$,
then $V$ is said to be {\em $A$-admissible} provided $V$ is an invariant subspace for each $S_A(t)$, $t\geq 0$, and the family of restrictions $\{S_{A}(t)_V:V\to V\}_{t\geq 0}$ ($S_A(t)_V \bar v:=S_A(t)\bar v$, $\bar v\in V$)
is a $C_0$ semigroup on $V$. Define the \emph{part of $A$ in the space $V$} as a linear operator $A_V:D(A_V)\to V$ given by $
D(A_V):=\{\bar v\in D(A)\cap V\,\mid\, A \bar v\in V\}$, $A_V \bar v:=A \bar v$ for $\bar v\in D(A_V)$.
In view of \cite[Ch. 4, Theorem 5.5]{Pazy}, if $V$ is $A$-admissible then $A_V$ is the generator of the $C_0$ semigroup $\{S_A(t)_V\}_{t\geq 0}$.\\
\indent Now let $\{A(t)\}_{t\geq 0}$ be a family of linear operators in $E$ satisfying the following conditions\\[2mm]
\noindent $(Hyp_1)$ \parbox[t]{140mm}{$\{A(t)\}_{t\geq 0}$ is
a {\em stable family} of infinitesimal generators of $C_0$ semigroups on $E$, i.e., there are $M\geq 1$ and $\omega\in\R$ such that
$$
\|S_{A(t_n)}(s_n)\ldots S_{A(t_1)}(s_1)\|_{{\cal L}(E,E)}\leq M e^{\omega(s_1+\ldots+s_n)},
$$
whenever $0\leq t_1 \leq \ldots \leq t_n$ and $s_1,\ldots, s_n\geq 0$, where $\{S_{A(t)}(s)\}_{s\geq 0}$ is the $C_0$ semigroup generated by $A(t)$;} \\[2mm]
\noindent $(Hyp_2)$ \parbox[t]{140mm}{$V$ is $A(t)$-admissible for each $t\geq 0$ and the family $\{ A_V (t)\}_{t\geq 0}$ is a stable family of generators of
$C_0$ semigroups on $V$ with constants $M_V\geq 1$ and $\omega_V \in\R$;} \\[1mm]
\noindent $(Hyp_3)$ \parbox[t]{140mm}{$V\subset D(A(t))$ and $A(t)\in{\cal L}(V,E)$, for any $t\geq 0$, and
the mapping $[0,+\infty) \ni t\mapsto A(t)\in {\cal L}(V,E)$ is continuous.}\\[1mm]
\noindent These are so called {\em hyperbolic conditions} and they determine a unique evolution system.
\begin{Prop}{\em (see \cite[Ch. 5, Theorem 3.1]{Pazy})}\label{17092008-1611}
Let $\{A(t)\}_{t\geq 0}$ be a family of linear operators in a Banach space $E$ satisfying $(Hyp_1)-(Hyp_3)$.
Then there exists a unique evolution system $\{R(t,s)\}_{t\geq s\geq 0}$ on $E$ with the following properties\\[1mm]
\indent {\em (i)} \parbox[t]{138mm}{$\|R(t,s)\| \leq M e^{\omega(t-s)}$ \quad for  $s\geq 0$;}\\[1mm]
\indent {\em (ii)} \parbox[t]{138mm}{$\left.\frac{\partial^+}{\partial t} R(t,s)v\right|_{t=s}= A(s)v$
\quad for \ $v\in V$, $s\geq 0$;}\\[1mm]
\indent {\em (iii)} \parbox[t]{138mm}{$\frac{\partial}{\partial s} R(t,s)v = -R(t,s)A(s)v$ \quad for \ $v\in V$, $0\leq s\leq t$.}
\end{Prop}
We shall consider parameterized evolution systems.
\begin{Prop}\label{09032010-1941}
Let $P$ be a metric space of parameters. Suppose that families $\{A^{(\mu)}(t)\}_{t\geq 0}$, $\mu\in P$,
satisfy conditions $(Hyp_1)-(Hyp_3)$ with constants $M, M_V, \omega, \omega_V$ independent of $\mu$. Let $R^{(\mu)} = \{ R^{(\mu)} (t,s) \}_{t\geq s\geq 0}$
be the corresponding evolution systems on $E$ determined by Proposition {\em \ref{17092008-1611}}. \\
\indent {\em (i)} \parbox[t]{140mm}{For any $\mu,\nu\in P$, $\bar v\in V$ and $t,s\geq 0$ with $t\geq s$,
$$
\|R^{(\nu)}\! (t,s)\bar v\!-\! R^{(\mu)}\!(t,s)\bar v\|\! \leq\!
M \! M_Ve^{(|\omega| +|\omega_V|)t}\|\bar v\|_{V} \!\int_{s}^{t}\!\!\! \|A^{(\nu)}(r)\!-\!A^{(\mu)}(r)\|_{{\cal L}(V,E)} \d r.
$$
}\\
\indent {\em (ii)} \parbox[t]{140mm}{If
\begin{equation}\label{13072009-1426}
\lim_{\nu\to \mu} \!\int_{0}^{T}\! \|A^{(\nu)}(r) - A^{(\mu)} (r)\|_{{\cal L}(V,E)} \d r=0, \   \mbox{ for any $\mu\in P$
and $T>0$,}
\end{equation}
then $\{ R^{(\mu)}\}_{\mu\in P}$ is a continuous family of evolution systems on $E$.}\\
\end{Prop}
{\bf Proof:} (i) We use the construction from \cite[Ch. 5, Theorem
3.1]{Pazy}. Fix $T>0$. Recall that, for any $\mu\in P$ and $\bar
u\in E$,
$$
R^{(\mu)}(t,s)\bar u = \lim_{n\to +\infty} R_{n}^{(\mu)} (t,s)\bar u \quad\mbox{ for \ } 0\leq s\leq t \leq T,
$$
where, for each $n\geq 1$, the operator $R_{n}^{(\mu)} (t,s):E\to
E$ is given by (\footnote{Here we use the convention that
$\prod_{j=m}^n T_j := T_n\circ T_{n-1}\circ\ldots  \circ T_m$, for
integers $m,n$ with $m<n$ and a sequence $T_m, T_{m+1}, \ldots, T_n$ of bounded operators on $E$.})
$$ R_{n}^{(\mu)}(t,s)\!:= \!\!\left\{
\begin{array}{ll}
\!\! S^{(\mu)}_j(t\!-\!s) &  \mbox{ if }  s,t\in [t_{j}^{n}, t_{j+1}^{n}], \, s \le t,\\
\!\! S^{(\mu)}_k(t\!-\!t_{k}^{n})\!\left(\prod\limits_{j=l+1}^{k-1} \!\!S^{(\mu)}_{j}(T/n)\! \right)\!\!S^{(\mu)}_{l}(t_{l+1}^{n}\!-\!s)&
\mbox{ if } s\!\in\! [t_{l}^{n}, t_{l+1}^{n}], \, t\!\in\! [t_{k}^{n}, t_{k+1}^{n}]\\
&\mbox{ and } k>l \geq 0,
\end{array} \right.
$$
with $t_{j}^{n}:=(j/n)T$, $S_{j}^{(\mu)}:=S_{A^{(\mu)}(t_{j}^n)}$, for $j=0,1,\ldots, n$.
Moreover $\{R_{n}^{(\mu)}(t,s)\}_{0\leq s\leq t\leq T}$, $\mu\in P$,  are evolution systems such that
\be\label{14012009-1920}
 \ \ \ \|R_{n}^{(\mu)}\!(t,s)\|_{{\cal L}(E,E)}
\! \leq\! Me^{\omega(t\!-\!s)}, \ \   R_n^{(\mu)}\!(t,s)V\!\subset \! V, \ \
\|R_{n}^{(\mu)}\!(t,s)\|_{{\cal L} (V,V)} \!\leq \! M_V e^{\omega_V
(t\!-\!s)},
\ee
whenever $0\leq s\leq t\leq T$ and, for any $\bar v\in V$,
\be\label{14012009-1921}
\frac{\partial}{\partial t} R_{n}^{(\mu)}(t,s)\bar v
=A_{n}^{(\mu)}(t)R_n^{(\mu)}(t,s)\bar v  \mbox{ for } t\not\in \{t_{0}^{n},t_{1}^{n},\ldots, t_{n}^{n}\}, \ s\leq t,
\ee
\be\label{14012009-1922} \frac{\partial}{\partial s}
R_{n}^{(\mu)}(t,s)\bar v  = - R_n^{(\mu)}(t,s) A_{n}^{(\mu)}(s)\bar v \mbox{ for \ }s\not\in
\{t_{0}^{n},t_{1}^{n},\ldots, t_{n}^{n}\}, \ s\leq t,
\ee
with $A_{n}^{(\mu)}(t):=A^{(\mu)}(t_{k}^{n})$ if $t_k^n \leq t <
t_{k+1}^{n}$ for $k=0,\ldots,n-1$ and
$A_n^{(\mu)}(T):=A^{(\mu)}(T)$.
$(Hyp_3)$ implies that, for any $\mu\in P$, one has
\be\label{15062010-2055}
\|A_{n}^{(\mu)}(t) -
A^{(\mu)}(t)\|_{{\cal L}(V,E)} \to 0 \mbox{ as  } n\to +\infty \mbox{ uniformly with respect to } t\in [0,T].
\ee
Fix any $\bar v\in V$, $\mu,\nu \in P$, $n\geq 1$ and
$s,t \in [0,T]$ with $s<t$ and define
$\phi:[s,t]\to E$ by $\phi(r):= R_n^{(\mu)}(t,r)R_n^{(\nu)}(r,s)\bar v$. In view of
(\ref{14012009-1920}), (\ref{14012009-1921}) and
(\ref{14012009-1922}), the map $\phi$ is differentiable on $[s,t]$
except the finite number of points and
\begin{eqnarray*}
R_n^{(\nu)}(t,s)\bar v- R_n^{(\mu)}(t,s)\bar v = \phi(t)-\phi(s)
=\int_{s}^{t} \phi'(r) \d r \\
= \int_{s}^{t}  R_n^{(\mu)}(t,r)(A_n^{(\nu)}(r))-A_n^{(\mu)}(r))R_n^{(\nu)}(r,s)\bar v\, \d r.
\end{eqnarray*}
This together with (\ref{14012009-1920}) yields
$$
\| R_n^{(\nu)}(t,s)\bar v- R_n^{(\mu)}(t,s)\bar v \|
\leq MM_Ve^{(|\omega| +|\omega_V|)t}\|\bar v\|_{V} \int_{s}^{t} \|A_n^{(\nu)}(r)-A_n^{(\mu)}(r)\|_{{\cal L}(V,E)} \d r.
$$
Passing to the limit with $n\to +\infty$, we get
$$
\|R^{(\nu)} (t,s)\bar v - R^{(\mu)}(t,s)\bar v\| \leq
MM_Ve^{(|\omega| +|\omega_V|)t}\|\bar v\|_{V} \int_{s}^{t} \|A^{(\nu)}(r)-A^{(\mu)}(r)\|_{{\cal L}(V,E)} \d r.
$$
\indent (ii) It follows immediately from (i) that $R^{(\nu)}(t,s)\bar v\to R^{(\mu)}(t,s)\bar v$ for any $\bar v\in V$,
as $\nu\to \mu$ and the convergence is uniform with respect to $s,t$ from $[0,T]$. To see it for an arbitrary
$\bar u\in E$, note that, for any $\bar v\in V$,
\begin{eqnarray*}
\|R^{(\nu)} (t,s)\bar u  - R^{(\mu)}(t,s)\bar u\| \leq \|R^{(\nu)} (t,s)\bar v  - R^{(\mu)}(t,s)\bar v\| + 2Me^{|\omega| T} \|\bar u - \bar v\|,
\end{eqnarray*}
which, in view of the density of $V$ in $E$, implies (ii).\hfill $\square$

Suppose that $\{A^{(\mu)}(t)\}_{t\geq 0}$, $\mu\in P$, are as in Proposition \ref{09032010-1941}.
Define $A^{(\mu, \lma)}(t): D(A^{(\mu, \lma)}(t)) \to E$, for $\lma>0$ and $t\geq 0$, by
$$
A^{(\mu,\lma)} (t):= A^{(\mu)}(t/\lma).
$$
Note that, for each $\mu\in P$ and $\lma>0$, the family $\{A^{(\mu,\lma)}(t) \}_{t\geq 0}$
also satisfies $(Hyp_1)-(Hyp_3)$ with the same constants (independent of $\mu$, $\lma$) as for $\{A^{(\mu)}(t) \}_{t\geq 0}$, $\mu\in P$. The following linear averaging principle holds.
\begin{Th}\label{02032010-0959}
Let $\{A^{(\mu)}(t) \}_{t\geq 0}$, $\mu\in P$, be as in
Proposition {\em \ref{09032010-1941}}. If, additionally, for each
$\mu\in P$, there is a generator $\widehat A^{(\mu)}$ of a $C_0$
semigroup $\{\widehat S^{(\mu)}(t):E\to E\}_{t\geq 0}$ such that
$V$ is $\widehat A^{(\mu)}$-admissible, $V\subset D(\widehat
A^{(\mu)})$ and
\be\label{02032010-0954} \ \ \ \lim_{T\to +\infty,\, \nu \to \mu} \frac{1}{T} \int_{0}^{T}\! \!
\| \!A^{(\nu)} (t+h) - \widehat A^{(\mu)}\|_{{\cal L}(V,E)} \d t  =0 \mbox{
uniformly with respect to } h\geq 0,
\ee
then, for any $\mu\in P$ and $\lma>0$,
\be\label{05032010-1223} \|R^{(\mu,\lma)}(t,s)\|
\leq M e^{\omega (t-s)} \mbox{ whenever } 0\leq s\leq t \ee
and, for any $t\geq 0$, $s\in [0,t]$, $\mu\in P$ and $\bar u\in E$,
$$
\lim\limits_{\lma \to 0^+,\ \bar v \to \bar u,\ \nu \to \mu} R^{(\nu, \lma)} (t,s)\bar v = \widehat S^{(\mu )}(t-s)\bar u
$$
where the convergence is uniform for $t$ and $s$ from bounded
intervals.
\end{Th}
\begin{Lem}\label{05032010-1004}
Under the assumption {\em (\ref{02032010-0954})}, for any $\mu \in P$ and $t,s\geq 0$ with $t\geq s$,
$$
\lim_{\lma\to 0^+, \, \nu\to\mu} \int_{s}^{t} \|A^{(\nu,\lma)} (r)-\widehat A^{(\mu)} \|_{{\cal L}(V,E)} \d r=0.
$$
\end{Lem}
{\bf Proof:} Take any $\mu\in P$ and let $(\lma_n)$ in
$(0,+\infty)$ and $(\mu_n)$ in $P$ be arbitrary sequences such
that $\lma_n\to 0$ and $\mu_n\to \mu$ in $P$. Let a sequence
$(k_n)$ of positive integers and $(T_n)$ in $(0,+\infty)$ be such
that $k_n\to +\infty$, $T_n\to +\infty$, $k_n \lma_n T_n\to t-s$
as $n\to +\infty$ and $k_n \lma_n T_n \leq t-s$ for each $n\geq
1$. Observe that
\begin{eqnarray*}
& & \int_{s}^{t} \|A^{(\mu_n,\lma_n)} (r) - \widehat A^{(\mu)}\| \d r = \int_{0}^{t-s}  \|A^{(\mu_n,\lma_n)} (s+\rho) - \widehat A^{(\mu)}\|\, \d\rho\\
& & = \int_{0}^{k_n\lma_n T_n}  \|A^{(\mu_n,\lma_n)} (s+\rho) - \widehat A^{(\mu)}\|\, \d\rho
+ \int_{k_n\lma_n T_n}^{t-s}  \|A^{(\mu_n,\lma_n)} (s+\rho) - \widehat A^{(\mu)}\|\, \d\rho.
\end{eqnarray*}
Clearly the second term tends to $0$ as $n\to +\infty$. Furthermore
\begin{eqnarray*}
& &\int_{0}^{k_n\lma_n T_n}\|A^{(\mu_n,\lma_n)} (s+\rho) - \widehat A^{(\mu)}\|\, \d\rho = \sum_{k=0}^{k_n-1} \int_{k\lma_n T_n}^{(k+1)\lma_n T_n}  \|A^{(\mu_n,\lma_n)} (s+\rho) - \widehat A^{(\mu)}\|\, \d\rho \\
& &\mbox{ }\mbox{ } = \lma_n T_n \sum_{k=0}^{k_n-1} \frac{1}{T_n} \int_{0}^{T_n}\|A^{(\mu_n)} (s/\lma_n+kT_n + \tau) -\widehat A^{(\mu)} \| \d \tau
\end{eqnarray*}
and, in view of (\ref{02032010-0954}), it also converges to $0$ as $n\to +\infty$.  \hfill $\square$

\noindent {\bf Proof of Theorem \ref{02032010-0959}:}
First, observe that the families of operators $\{ A^{(\mu,\lma)}(t)\}_{t\geq 0}$, $\mu\in P$, $\lma>0$, also satisfy
$(Hyp_1)-(Hyp_3)$ and, by the very construction, the corresponding evolution system $\{R^{(\mu,\lma)}(t)\}_{t\geq 0}$
admits the same growth condition as $\{R^{(\mu)}(t,s)\}_{t\geq s\geq 0}$, i.e., (\ref{05032010-1223}) holds.
Next we intend to prove that, for any $t, s\geq 0$ with $t\geq s$, $\mu\in P$ and $\bar u\in E$, the limit
$\lim_{\lma\to 0^+, \nu\to \mu, \bar v\to \bar u} R^{(\nu, \lma)} (t,s)\bar v$ exists.  To this end fix any $T>0$ and note
that it follows from Proposition
\ref{09032010-1941} (i) that, for $C:= MM_Ve^{(|\omega| +|\omega_V|)T}$, any $\mu_0, \mu_1,\mu_2\in P$, $\lma_1, \lma_2>0$, $\bar v\in V$
and $t,s\in [0,T]$ with $t\geq s$, one has
\begin{eqnarray*}
& & \|R^{(\mu_1,\lma_1)}(t,s)\bar v - R^{(\mu_2,\lma_2)}(t,s)\bar v\|
\leq C \|\bar v\|_{V} \int_{0}^{T} \|A^{(\mu_1,\lma_1)}(r)-A^{(\mu_2,\lma_2)}(r)\|_{{\cal L}(V,E)} \d r\\
& & \ \ \ \ \ \ \ \ \
\leq C \|\bar v\|_{V} \left( \int_{0}^{T} \|A^{(\mu_1,\lma_1)}(r)-\widehat A^{(\mu_0)}\|_{{\cal L}(V,E)} \d r + \right.\\
& & \ \ \ \ \ \ \ \ \ \ \ \ \ \ \ \ \ \  \ \ \ \ \ \ \ \ \ \left. +\int_{0}^{T} \|\widehat A^{(\mu_0)}- A^{(\mu_2, \lma_2)}(r)\|_{{\cal L}(V,E)} \d r \right).
\end{eqnarray*}
Now it is immediate, in view of Lemma \ref{05032010-1004}, that the
limit $\lim_{\lma\to 0^+,\, \nu\to \mu} R^{(\mu,\lma)}(t,s)\bar v$ exists
and is uniform with respect to $t,s$ from bounded intervals. If one takes any $\bar u \in E$,
then the existence of the limit $\lim_{\lma\to 0^+,\, \nu\to \mu} R^{(\nu,\lma)}(t,s)\bar u$ and its uniformness with respect to $t,s\in [0,T]$ come from
the inequality
$$
\|R^{(\mu_1,\lma_1)} (t,s) \bar u  - R^{(\mu_2,\lma_2)} (t,s) \bar u\| \leq
\|R^{(\mu_1,\lma_1)} (t,s) \bar w  - R^{(\mu_2,\lma_2)} (t,s) \bar w\| + 2Me^{|\omega| T} \|\bar u - \bar w\|
$$
holding for any $\mu_1,\mu_2\in P$, $\lma_1,\lma_2>0$ and arbitrary $\bar w\in V$.
Finally, by the growth condition (\ref{05032010-1223}), for $\mu_1,\mu_2\in P$, $\lma_1,\lma_2>0$ and $\bar v_1, \bar v_2 \in E$,
\begin{eqnarray*}
& & \| R^{(\mu_1, \lma_1)} (t,s) \bar v_1 - R^{(\mu_2, \lma_2)} (t,s)\bar v_2 \| \leq
\|R^{(\mu_1, \lma_1)} (t,s) \bar v_1  - R^{(\mu_1, \lma_1)} (t,s)\bar u \| \\
& & \ \ \ \ + \|R^{(\mu_1, \lma_1)} (t,s)\bar u - R^{(\mu_2, \lma_2)} \bar u \| + \| R^{(\mu_2, \lma_2)} \bar u  - R^{(\mu_2, \lma_2)} (t,s)\bar v_2 \|\\
& & \ \ \ \ \leq M e^{|\omega|T} \| \bar v_1 - \bar u\| + \|R^{(\mu_1, \lma_1)} (t,s)\bar u - R^{(\mu_2, \lma_2)} \bar u \| +
M e^{|\omega|T} \| \bar v_2 - \bar u\|,
\end{eqnarray*}
which implies that $\lim_{\lma\to 0^+,\, \nu\to \mu,\, \bar v\to \bar u} R^{(\nu,\lma)} (t,s)\bar v$ exists and is uniform with respect to $t,s$ from bounded intervals.\\
\indent Let operators $\widehat R^{(\mu)} (t,s): E\to E$, $t,s\geq 0$ with $t\geq s$, $\mu\in P$, be defined by
$$
\widehat R^{(\mu)}(t,s)\bar u := \lim_{\lma\to 0^+} R^{(\mu,\lma)} (t,s)\bar u, \mbox{ for any } \bar u \in E.
$$
To complete the proof we need to show that $\widehat R^{(\mu)}(t,s) = \widehat S^{(\mu)}(t-
s)$ for any $t, s\geq 0$ with $t\geq s$ and $\mu\in P$.
For fixed $\mu\in P$, $\lma>0$, $n\geq 1$,  $\bar v\in V$, $T>0$ and $t,s\in [0,T]$ with $t\geq s$,
define $\varphi_\lma :[s,t]\to E$ by
$$
\varphi_{\lma} (r):=\widehat S^{(\mu)} (t-r) R_{n}^{(\mu,\lma)}(r,s)\bar v
$$
where $R_{n}^{(\mu,\lma)}$ are approximating evolution systems from the construction of $R^{(\mu,\lma)}$
(see the proof of Proposition \ref{09032010-1941}). Since $R_{n}^{(\mu,\lma)}(r,s)\bar v\in V\subset D(\widehat A^{(\mu)})$,
one has, for a.e. $r\in [s,t]$,
\begin{eqnarray*}
\varphi'_{\lma}(r)& =& \frac{\part}{\part r} \left( \widehat S^{(\mu)}(t-r)R_{n}^{(\mu,\lma)}(r,s)\bar v \right)
= [-\widehat A^{(\mu)} \widehat S^{(\mu)} (t-r)] R_n^{(\mu,\lma)}(r,s)\bar v  \\
& & +  \widehat S^{(\mu)}(t-r) [ A_{n}^{(\mu,\lma)}(r) R_{n}^{(\mu,\lma)}(t,r)\bar v]\\
& = & \widehat S^{(\mu)} (t-r) ( A_{n}^{(\mu,\lma)}(r) -  \widehat A^{(\mu)} ) R_{n}^{(\mu,\lma)}(r,s) \bar v
\end{eqnarray*}
(cf. (\ref{14012009-1921})). Hence, using the estimates for
$R_n^{(\mu,\lma)}$ as those in (\ref{14012009-1920}), one gets
\begin{eqnarray*}
& & \| R_{n}^{(\mu,\lma)}(t,s)\bar v-\widehat S^{(\mu)}(t-s)\bar v \| = \|\varphi_\lma (t)-\varphi_\lma  (s)\| \leq
 \int_{s}^{t} \|\varphi'_\lma (r)\|\, \d r \\
 & & \ \ \ \ \leq
 \int_{s}^{t}  \|\widehat S^{(\mu)} (t-r)\|_{{\cal L}(E,E)} \| A_{n}^{(\mu,\lma)}(r) -  \widehat A^{(\mu)} \|_{{\cal L}(V,E)} \|R_{n}^{(\mu,\lma)}(r,s)\|_{{\cal L}(V,V)} \|\bar v\|_{V} \d r \\
 & &  \leq \widehat M e^{|\widehat \omega| t} M_Ve^{|\omega_V|t} \|v\|_{V}
\int_{s}^{t} \|A_{n}^{(\mu,\lma)} (r)-  \widehat A^{(\mu)} \|_{{\cal L}(V,E)} \, \d r
 \end{eqnarray*}
 where $\widehat M\geq 0$ and $\widehat \omega\in\R$ are such that $\|\widehat S^{(\mu)}(\tau)\|_{{\cal L}(E,E)}\leq \widehat M e^{\widehat \omega \tau}$ for any $\tau\geq 0$.
Letting $n\to +\infty$, one obtains
$$\|\widehat S^{(\mu)}(t-s)\bar v - R^{(\mu,\lma)}(t,s)\bar v\| \leq
\widehat M M_Ve^{(|\widehat \omega|+ |\omega_V|)t} \|v\|_{V}
\!\int_{s}^{t}\! \|A^{(\mu,\lma)}(r) -  \widehat A^{(\mu)} \|_{{\cal L}(V,E)} \, \d r
$$
(cf. (\ref{15062010-2055})). In view of Lemma \ref{05032010-1004}, if we pass to the limit with $\lma\to 0^+$ then
$\widehat S^{(\mu)}(t-s)\bar v = \widehat R^{(\mu)}(t,s)\bar v$. Since $\bar v\in V$
was arbitrary and $V$ is dense in $E$, we find that $\widehat S^{(\mu)}(t-s) = \widehat R^{(\mu)}(t,s)$. \hfill $\square$

\section{Example}
Consider the following system of equations
$$
\frac{\part u}{\part t}(x,t) = \sum_{j=1}^{N} a_j (x,t) \frac{\part u}{\part x_j}(x,t) + b(x,t)u(x,t) + f(x,t), \ \ \ x\in\R^N, t>0
$$
with continuously differentiable $a_j: \R^N \times [0,+\infty) \to {\cal M}$, $j=1, \ldots , N$ ($N\geq 1$), and continuous $b:\R^N\times [0,+\infty)
\to {\cal M}$ with ${\cal M}$ being the space of all square matrices of order $M\geq 1$ with the usual maximum norm denoted by $|\cdot|$
and $u:\R^N\times [0,+\infty)\to \R^M$.
We assume that $a_j(x,t)$ is
symmetric for any $j=1,\ldots, N$ and $(x,t)\in\R^N\times [0,+\infty)$.
By $B^{0}(\R^N, {\cal M})$ denote the space of
all bounded continuous $\varphi:\R^N \to  {\cal M} $ with the norm $\| \varphi\|_{B^0 (\R^N, {\cal M})}:=
\sup_{x\in\R^N} |\varphi(x)|.$ And let $B^{1}(\R^N, {\cal M})$ denote the space of all continuously differentiable functions
$\varphi:\R^N \to  {\cal M} $ such that both $\varphi$ and its derivative are bounded
and equip it with the norm $\| \varphi\|_{B^1(\R^N, {\cal M})}
:= \sup_{x\in\R^N} \left(|\varphi(x)|+\|\varphi'(x)\|_{{\cal L} (\R^N, {\cal M})}\right)$.\\
\indent We shall also assume that the maps $t\mapsto a_j (\cdot,
t) \in B^{1}(\R^N, {\cal M})$, $j=1,\ldots, N$, and $t\mapsto b
(\cdot, t) \in B^{0}(\R^N, {\cal M})$ are well defined, continuous
and bounded. By use of \cite[Section 4.6]{Tanabe} we conclude that
the family of operators $\{ A(t) \}_{t\geq 0}$ in $E:= L^2 (\R^N,\R^M)$ given by
$$
D(A(t)):=\left\{ u\in E \mid \frac{\partial u}{\partial x_j} \mbox{ exists for each }j=1,\ldots, N, \sum_{j=1}^{N} a_j (\cdot,t) \frac{\part u}{\part x_j} \in E      \right\},
$$
where weak derivatives are considered, and
$$
[A(t)u](x)\!:= \!\sum_{j=1}^{N}\! a_j (x,t) \frac{\part u}{\part x_j}(x) + b(x,t)u(x), \mbox{ for each } u\in D(A(t)), \mbox{ a.e. } x\in \R^N,
$$
satisfies $(Hyp_1)-(Hyp_3)$  with $V:=H^1(\R^N,\R^M)$, $M=M_{V}=1$
and some constants $\omega, \omega_V$. Consequently, due to Proposition \ref{17092008-1611}, the family $\{A(t)\}_{t\geq 0}$ determines
a unique evolution system $\{R(t,s)\}_{t\geq s\geq 0}$ on $E$.
Therefore, if we assume that $f:\R^N\times [0,+\infty)\to\R^M$ is such that the map
$[0,+\infty)\ni t\mapsto f(\cdot, t) \in L^2(\R^N,\R^M)$ is well defined and continuous, then, for $\lma>0$, we may consider
\be\label{30042010-0835}
\dot u(t) = A(t/\lma)u(t)+F(t/\lma)
\ee
with $F:[0,+\infty) \to E$ given by $F(t):=f(\cdot ,t)$ for each $t\geq 0$.\\
\indent Now suppose that there are $\widehat a_j\in B^{1}(\R^N, {\cal M})$, $j:=1,\ldots,N$, with symmetric values
and $\widehat b\in B^{1}(\R^N, {\cal M})$ such that
\be\label{09032010-2255}
\lim\limits_{T\to + \infty} \frac{1}{T}\int_{0}^{T} \|a_j (\cdot,\tau + h) - \widehat a_j\|_{B^0 (\R^N,{\cal M})}
\, \d\tau = 0
\ee
and
\be\label{09032010-2256}
\lim\limits_{T\to + \infty} \frac{1}{T}\int_{0}^{T}\|b (\cdot,\tau + h )-\widehat b\|_{B^0(\R^N,{\cal M})}   \, \d\tau = 0
\ee
uniformly with respect  $h>0$. We claim that
\be\label{10032010-1341}
\lim_{T\to +\infty} \frac{1}{T}\int_{0}^{T}\|A(\tau + h)-\widehat A \|_{{\cal L}(V,E)}\, \d\tau = 0
\ee
uniformly with respect to $h>0$, where $\widehat A:D(\widehat A) \to E$ is defined by
$$
D(\widehat A):=\left\{ u\in E \mid \frac{\partial u}{\partial x_j} \mbox{ exists for each }j=1,\ldots, N, \sum_{j=1}^{N} \widehat a_j \frac{\part u}{\part x_j} \in E      \right\},
$$
$$
\widehat A u:=\sum_{j=1}^{N} \widehat a_j \frac{\part u}{\part x_j} + \widehat b u.
$$
Indeed, for any $u\in V$ and $t\geq 0$,
\begin{eqnarray*}
& & \|A(t)u-\widehat A u\|_{E}  \leq  \displaystyle{ \left(
\sum_{j=1}^{N} \int_{\R^N} |a_j (x,t) \frac{\part u}{\part x_j}(x) -
\widehat a_j (x) \frac{\part u}{\part x_j}(x)  |^2 dx \right)^{1/2}} \\
& & \ \ \ \ \ \  \displaystyle{+ \left( \int_{\R^N} |b(x,t)u(x)-\widehat b(x)u(x)|^2 dx \right)^{1/2}} \\
& & \leq
\left(\displaystyle{\sum_{j=1}^{N} \sup_{x\in \R^N} |a_j (x ,t) - \widehat a_j (x) |  +
\sup_{x\in\R^N} |b(x,t)-\widehat b(x)| }\right)\|u\|_{V}
\end{eqnarray*}
and this yields, for any $h>0$,
\begin{eqnarray*}
\displaystyle{\frac{1}{T}\int_{0}^{T}\|A(\tau+h)-\widehat A \|_{{\cal L}(V,E)}\, \d\tau \leq
\sum_{j=1}^{N}  \frac{1}{T}\int_{0}^{T} \|a_j (\cdot , \tau +h)- \widehat a_j \|_{B^0(\R^N,{\cal M})} \, \d\tau} \\
\ \ \ \ \ \ \ \ + \displaystyle{\frac{1}{T}\int_{0}^{T} \|b (\cdot,\tau+h) - \widehat b\|_{B^0 (\R^N,{\cal M})}  \, \d\tau,}
\end{eqnarray*}
which, by use of (\ref{09032010-2255}) and (\ref{09032010-2256}),
gives (\ref{10032010-1341}) uniformly with respect to $h>0$.
Consequently, if $\{R^{(\lma)}(t,s) \}_{t\geq s\geq 0}$, $\lma>0$,
denote the evolution system generated by $\{ A(t/\lma)\}_{t\geq
0}$ and $\{\widehat S(t)\}_{t\geq 0}$ the semigroup generated by
$\widehat A$, then, in view of Theorem \ref{02032010-0959}, we get
that, for any $\bar u\in E$, $\mu\in P$ and $t,s\geq 0$ with
$t\geq s$,
$$
\lim_{\lma\to 0^+, \bar v\to \bar u} R^{(\lma)}(t,s)\bar v = \widehat S(t-s)\bar u
$$
uniformly with respect to $t,s$ from bounded intervals.\\
\indent At last suppose that
$$
\int_{\R^N} |f(x+y,t)-f(x,t)|^2 \d x \to 0 \mbox{ as } y\to 0, \mbox{ uniformly in } t\geq 0,
$$
$$
\int_{\R^N \setminus B(0,R)} |f(x,t)|^{2} \d x \to 0 \mbox{ as } R\to 0^+, \mbox{ uniformly in } t\geq 0,
$$
and that there is $\widehat f\in L^{2}(\R^N,\R^M)$ such that
$$
\lim_{T\to +\infty} \int_{\R^N} \left| \frac{1}{T} \int_{0}^{T} f(x,t+h) \d t -\widehat f(x)\right|^2 \d x = 0, \mbox{ uniformly with respect to } h>0.
$$
Then the set $F([0,+\infty))$ is relatively compact in $L^2 (\R^N, \R^M)$ and
$$
\lim_{T\to +\infty} \frac{1}{T} \int_{0}^{T} F(t+h) \d t = \widehat f \mbox{ in } L^2 (\R^N,\R^M) \mbox{ uniformly with respect to } h>0.
$$
Thus, by use of Theorem \ref{25112009-1155} we conclude that mild solutions of
(\ref{30042010-0835}) converge to mild solutions of the averaged equation $\dot u = \widehat A u+ \widehat f$, i.e.,
the system
$$
\frac{\part u}{\part t}(x,t) = \sum_{j=1}^{N} \widehat a_j (x) \frac{\part u}{\part x_j}(x,t) + \widehat b(x)u(x,t) + \widehat f(x), \ \ \ x\in\R^N, t>0.
$$

\end{document}